\documentclass[amstex,12pt, amssymb]{article}

\usepackage{mathtext}
\usepackage[cp1251]{inputenc}
\usepackage[T2A]{fontenc}
\usepackage[dvips]{graphicx}
\usepackage{amsmath}
\usepackage{amssymb}
\usepackage{amsxtra}
\usepackage{latexsym}
\usepackage{ifthen}

\newcounter{lemma}[section]

\newcounter{corollary}[section]

\newcounter{remark}[section]

\newcounter{theorem}[section]

\newcounter{proposition}[section]

\numberwithin{equation}{section}

\begin{document}

\markboth{\centerline{E. SEVOST'YANOV}} {\centerline{THE MINIOWITZ
AND VUORINEN THEOREMS \ldots }}

\author{E. SEVOST'YANOV\\}

\title{{\bf THE MINIOWITZ AND VUORINEN THEOREMS \\ FOR SOME CLASS OF MAPPINGS WITH
NON-BOUNDED CHARACTERISTICS }}

\author{{\bf E. Sevost'yanov}\\}
\date{\today \hskip 4mm }
\maketitle

\abstract The paper is devoted to studying classes of mappings with
unbounded characteristics of quasiconformality. Namely, we prove
that the normal families of $Q$-mappings have the logarithmic order
of growth in a neighborhood of a point. Moreover, we establish
sufficient conditions for $Q$ provided normality of mappings
$f:D\rightarrow \overline{{\Bbb R}^n},$ $n\ge 2,$ omitting points of
a set $E_f$ obeying $c(E_f)\ge \delta,$ $\delta>0,$ where $c(\cdot)$
is an appropriate set function.
\endabstract


\section{Introduction}\label{sect1}

The paper is devoted to studying ring $Q$-map\-pings, which
generalize $Q$-map\-pings, introduced in \cite{MRSY} (cf. \cite{RSY}
and \cite{GRSY} for homeomorphisms). The class of ring $Q$-mappings
contains many known classes of mappings such as analytic functions,
quasiconformal mappings, mappings with bounded and finite
distortion, etc.

Throughout the paper, $m$ denotes the Lebesgue measure in ${\Bbb
R}^n,$ $n\ge 2,$ $\overline{{\Bbb R}^n}={\Bbb R}^n\cup\{\infty\}$ is
the one--point compactification of ${\Bbb R}^n,$ and $M(\Gamma)$
stands for the conformal modulus of families of curves $\gamma$ in
${\Bbb R}^n$ (see e.g. \cite[6.I]{Va$_1$}). For $\overline{{{\Bbb
R}}^n},$ we use the {\it spherical (chordal)} distance
$h(x,y)=|\pi(x)-\pi(y)|$, where $\pi$ is a stereographical
projection of $\overline{{{\Bbb R}}^n}$ onto the sphere
$S^n(\frac{1}{2}e_{n+1},\frac{1}{2})$ in ${{\Bbb R}}^{n+1}$:
$$h(x,\infty)=\frac{1}{\sqrt{1+{|x|}^2}}, \ \
h(x,y)=\frac{|x-y|}{\sqrt{1+{|x|}^2} \sqrt{1+{|y|}^2}}\,, \quad x\ne
\infty\ne y\,.$$  Given a domain $D$ and two sets $E$ and $F$ in
${\overline{{\Bbb R}^n}},$ $n\ge 2,$ $\Gamma (E,F,D)$ denotes the
family of all paths $\gamma:[a,b]\rightarrow {\overline{{\Bbb
R}^n}}$ which join $E$ and $F$ in $D$, i.e., $\gamma(a)\in E,$
$\gamma(b)\in F$ and $\gamma(t)\in D$ for $a<t<b.$ We write
$\Gamma(E,F)= \Gamma(E,F,{\overline{{\Bbb R}^n}})$ when
$D={\overline{{\Bbb R}^n}}.$ Let $r_0=\,{\rm dist\,} (x_0\,,\partial
D)$ and $Q:D\rightarrow [0\,,\infty]$ be a measurable function. Set
$$R(r_1,r_2,x_0)= \{ x\,\in{\Bbb R}^n : r_1<|x-x_0|<r_2\}\,, $$
$$S_{\,i}=S(x_0,r_i)= \{ x\,\in\,{\Bbb R}^n: |x-x_0|=r_i\}\,,\quad i=1,2. $$
A mapping $f:D\rightarrow \overline{{\Bbb R}^n}$ is said to be a
{\it ring $Q$-map\-ping at a point $x_0\,\in\,D$} if the inequality
\begin{equation}\label{eq2} M\left(f\left(\Gamma\left(S_1, S_2,
R\right)\right)\right)\le \int\limits_{R} Q(x)\cdot \eta^n(|x-x_0|)\
dm(x) \end{equation}
holds for any $R=R(r_1,r_2, x_0),$\, $0<r_1<r_2< r_0,$ and every
measurable function $\eta: (r_1,r_2)\rightarrow [0,\infty ]$ with
\begin{equation}\label{eq2A}
\int\limits_{r_1}^{r_2}\eta(r) dr \ge 1\,.
\end{equation}
Note that for $Q(x)\le K= const,$ the class of $Q$-mappings contains
quasiregular mappings, whereas $1$-mappings have within themselves
conformal mappings (see \cite{Re$_2$}, \cite{Ri} and
\cite[Theorem~1]{Pol}).

\medskip
In \cite{Sev$_1$}, we establish conditions on $Q$ ensuring the
equicontinuity and normality of families of mappings satisfying
(\ref{eq2}) and omitting points of a fixed set of positive capacity.
These conditions are sufficient for equicontinuity but not
necessary. In particular, a family of mappings $f:D\rightarrow
\overline{{\Bbb R}^n}\setminus E$ obeying (\ref{eq2}) is normal
whenever ${\rm cap\,}E>0$ and $Q\in FMO$ (see
\cite[Theorem~5.1]{Sev$_1$}).

\medskip The paper consists of two main parts devoted to
studying the equicontinuity property for the mappings satisfying the
relations (\ref{eq2})--(\ref{eq2A}). In the first part we establish
necessary and sufficient conditions providing the equicontinuity for
a family of ring $Q$-mappings, imposing an appropriate condition on
$Q,$ e.g. $FMO$-condition. Note that our results generalize the
known results by Miniowitz for quasiregular mappings, see
\cite[Theorem~1]{M}. For the plane quasiregular mappings we refer to
\cite[Theorem~4.3.II]{LV} (cf. \cite{Zal}). The main result in this
part is the following

\medskip
\begin{theorem}{}\label{th3}{\,\sl A family of all discrete open ring
$Q$-map\-pings $f:D\rightarrow \overline{{\Bbb R}^n}$ at a point
$x_0\in D$ with $Q\in FMO(x_0)$ is equicontinuous at $x_0$ if and
only if there exist $p=p(n , Q)>0,$ $C_n>0$ and $\varepsilon_0\in
(0, {\rm\,dist}(x_0, \partial D))$ such that
\begin{equation}\label{eq1}
h(f(x), f(x_0))\le C_n
\left\{{\frac{1}{\log{\frac{1}{|x-x_0|}}}}\right\}^{p}
\end{equation} for every $x\in B(x_0, \varepsilon_0),$
where $B(x_0,\varepsilon)=\{x\in {\Bbb R}^n:
|x-x_0|<\varepsilon\}.$}
\end{theorem}

\medskip
Following \cite{IR}, we say that a function $\varphi:D\rightarrow
{\Bbb R} $ has {\it finite mean oscillation} ($FMO$) at a point $x_0
\in {D} $ if
$$\limsup\limits_{\varepsilon\rightarrow 0}\frac{1}{\Omega_n\varepsilon^n}
\int\limits_{B(x_0 ,\varepsilon)}
|\varphi(x)-\widetilde{\varphi_{\varepsilon}}|dm(x)< \infty\,,
$$
where
$$\widetilde{\varphi_{\varepsilon}}=
\frac{1}{\Omega_n\varepsilon^n}\int\limits_{B( x_0 ,\varepsilon)}
\varphi(x)\, dm(x)$$
%
is the mean of the function $\varphi(x)$ over the ball
$B(x_0,\varepsilon).$ Note that $FMO$ is not $BMO_{loc}$
(\cite[p.~211]{MRSY}). It is well known \cite{JN} that $L^{\,\infty
}(D)\subset BMO(D)\subset L^p_{loc}(D)$ for all $1\le p<\infty,$ but
$FMO(D)\not\subset L_{loc}^p(D)$ for any $p>1$ (see \cite{MRSY}).

\medskip
In the second part of the paper we establish sufficient conditions
of the equicontinuity for mappings $f:D\rightarrow \overline{{\Bbb
R}^n}\setminus E_f$ omitting the points of a set $E_f,$ where $E_f$
is a set of positive capacity depending on $f.$ Note that, in the
general case, the condition ${\rm cap}\,E_f>0$ does not imply the
equicontinuity (normality) of a family of $K$--qua\-si\-re\-gu\-lar
mappings $f:D\rightarrow \overline{{\Bbb R}^n}\setminus E_f.$ For
example, the family $f_m(z)=z^m,$ $z\in {\Bbb C},$ $D:=B(0,
2)=\{z\in {\Bbb C}: |z|<2\}$ for $n=2$ is not normal. It is known
that, if a set $E=E_f$ does not depend on $f,$ the corresponding
family of mappings is equicontinuous (see
\cite[Corollary~2.7.III]{Ri} for bounded $Q$ and
\cite[Theorems~5.1--5.2]{Sev$_1$} for more general $Q$). However, in
the present paper we do not require an existence of such
$"$general$"$\,\, $E.$ Following \cite{Vu$_1$}, we restrict the
family $\{E_f\}$ involving a set function $c(\cdot)$ (see the
relation (\ref{eq25})). The restrictions of such type have been
earlier applied for bounded $Q$ (cf. \cite{Vu$_1$}). The following
statement is one of the main results of the paper.

\medskip
\begin{theorem}{}\label{th1}
{\sl\, Let ${\frak F}_{Q, \delta}$ be a family of open discrete ring
$Q$-map\-pings $f:D\rightarrow\overline{{\Bbb R}^n}\setminus E_f$ at
$x_0\in D$ with $c(E_f)\ge \delta>0,$ where $E_f$ is compact. Assume
that either of the following conditions holds: 1) $Q\in FMO(x_0);$
2) $q_{x_0}(r)\le C\cdot (\log\frac{1}{r})^{n-1}$ where $C>0$ is a
positive constant and $q_{x_0}(r)$ is the integral mean of $Q(x)$
over the sphere $S(x_0, r);$ 3)
\begin{equation}\label{eq12}
\int\limits_{0}^{\varepsilon_0}
\frac{dt}{tq_{x_0}^{\,\frac{1}{n-1}}(t)}=\infty
\end{equation} for some $\varepsilon_0=\varepsilon_0(x_0).$
Then the family ${\frak F}_{Q, \delta}$ is equicontinuous at $x_0.$}
\end{theorem}

\medskip
\section{Preliminaries}\label{sect2}
\setcounter{equation}{0}

Recall the needed definitions and notions (see also \cite{Sev$_1$}).
Let $D$ be a domain in ${\Bbb R}^n,$ $n\ge 2.$ A mapping
$f:D\rightarrow {\Bbb R}^n$ is said to be {\it discrete} if the
preimage $f^{-1}\left(y\right)$ of every point $y\in{\Bbb R}^n$
consists of isolated points, and {\it open} if the image of every
open set $U\subset D$ is open in ${\Bbb R}^n.$

If $U$ is an open set in ${\Bbb R}^n,$ we denote by $C_0(U)$ the
class of all continuous functions in $U$ whose support is a compact
subset of $U.$ Following \cite{MRV$_1$} and \cite{Ri}, a {\it
condenser} is a pair $E=(A, C),$ where $A\subset {\Bbb R}^n$ is open
and $C$ is a non--empty compact set contained in $A.$ Define the
quantity
$$
{\rm cap}\,E={\rm cap}\,\left(A,\,C\right)=\inf \limits_{u\in
W_0\left(E\right) }\,\,\int\limits_{A}\,|\nabla u|^n dm(x)
$$
which is called the {\it capacity} of the condenser $E,$ where
$|\nabla u|={\left(\sum\limits_{i=1}^n\,{\left(\partial_i
u\right)}^2 \right)}^{1/2},$ $W_0(E)=W_0(A,\,C)$ is a family of
nonnegative functions $u: A\rightarrow {\Bbb R}^1$ such that
(1)\quad $u\in C_0(A),$ \quad (2)\quad $u(x)\ge 1$ for $x\in C$ and
(3)\quad $u$ is $ACL$ (absolutely continuous on lines, see
\cite{Ri}).

\medskip
Let $f:D \rightarrow {\Bbb R}^n$ be a discrete open mapping, $\beta:
[a,\,b)\rightarrow {\Bbb R}^n$ be a curve and $x\in
f^{\,-1}(\beta(a)).$ A curve $\alpha: [a,\,c)\rightarrow D$ is
called a {\it maximal $f$--lifting} of $\beta$ starting at $x$ if
$(1)\quad \alpha(a)=x;$ $(2)\quad f\circ\alpha=\beta|_{[a,\,c)};$
$(3)$\quad there exists no path $\alpha^{\,\prime}:
[a,\,c^{\prime})\rightarrow D$ such that
$\alpha=\alpha^{\prime}|_{[a,\,c)}$ and $f\circ
\alpha^{\,\prime}=\beta|_{[a,\,c^{\prime})}$ whenever
$c<c^{\prime}\le b.$ For a discrete open mapping $f,$ every curve
$\beta$ with $x\in f^{\,-1}\left(\beta(a)\right)$ has  a maximal
$f$--lifting starting at a point $x$ (see
\cite[Corollary~3.3.II]{Ri}, see also \cite[Lemma~3.12]{MRV$_3$}).

\medskip
The following statement was earlier formulated and exploited in more
special cases (see \cite[Lemmas~3.2--3.3]{Sev$_1$}).

\medskip
\begin{lemma}{}\label{lem4}{\,\sl
Let $f:D\rightarrow\overline{{\Bbb R}^n},$ $n\ge 2,$ be an open
discrete ring $Q$-map\-ping at a point $x_0\in D.$ Suppose that
there exist numbers $\varepsilon_0 \in (0,\, {\rm dist}\,(x_0,
\partial D)),$ $\varepsilon_0^{\,\prime}\in (0, \varepsilon_0)$ and
a family of nonnegative Lebesgue measurable functions
$\{\psi_{\varepsilon}(t)\},$ $\psi_{\varepsilon}:(\varepsilon,
\varepsilon_0)\rightarrow [0, \infty],$ $\varepsilon\in\left(0,
\varepsilon_0^{\,\prime}\right),$ such that
\begin{equation} \label{eq3.7B}
\int\limits_{\varepsilon<|x-x_0|<\varepsilon_0}Q(x)\cdot\psi_{\varepsilon}^n(|x-x_0|)
\ dm(x)\le F(\varepsilon,
\varepsilon_0)\qquad\forall\,\,\varepsilon\in(0,
\varepsilon_0^{\,\prime})\,,
\end{equation}
where $F(\varepsilon, \varepsilon_0)$ is a given function and
\begin{equation}\label{eq3} 0<I(\varepsilon, \varepsilon_0):=
\int\limits_{\varepsilon}^{\varepsilon_0}\psi_{\varepsilon}(t)dt <
\infty\qquad\forall\,\,\varepsilon\in(0,
\varepsilon_0^{\,\prime})\,.\end{equation}
Then
\begin{equation}\label{eq3B}
{\rm cap}\,f(E)\le F(\varepsilon,\varepsilon_0)/ I^{n}(\varepsilon,
\varepsilon_0)\qquad
\forall\,\,\varepsilon\in\left(0,\,\varepsilon_0^{\,\prime}\right)\,,
\end{equation}
where $E=\left(A,\,C\right),$ $A=B\left(x_0, r_0\right),$
$C=\overline{B(x_0, \varepsilon)},$ $r_0={\rm
dist}\left(x_0,\,\partial D\right),$ $A:={\Bbb R}^n$ whenever
$D:={\Bbb R}^n.$
}
\end{lemma}

\medskip
{\it Proof.} Consider a condenser $E=(A,\,C),$ where $A$ and $C$ are
defined as above. Obviously the pair $f(E)=(f(A), f(C))$ is also
condenser. If ${\rm cap}\,f(E)=0,$ (\ref{eq3B}) is trivial. Thus,
assume that ${\rm cap}\,f(E)\ne 0.$ One can assume that
$\infty\not\in f(A).$

\medskip
Let $\Gamma_E$ be a family of curves $\gamma:[a,\,b)\rightarrow A$
such that $\gamma(a)\in C$ and $|\gamma|\cap\left(A\setminus
F\right)\ne\varnothing$ for every compact set $F\subset A$ where
$|\gamma|=\{x\in {\Bbb R}^n: \exists\,t\in [a, b): \gamma(t)=x \}$
is a locus of $\gamma.$ It is known that ${\rm cap}\,E=M(\Gamma_E)$
(see, e.g. \cite[Proposition~10.2.II]{Ri}).

\medskip
Consider a family $\Gamma_{f(E)}$ for the condenser $f(E).$ Note
that every $\gamma\in\Gamma_{f(E)}$ has a maximal $f$--lif\-ting in
$A$ starting in $C$ (see \cite[Corollary~3.3.II]{Ri}). Let
$\Gamma^{*}$ be a family of all maximal $f$--lif\-tings of
$\Gamma_{f(E)}$ starting in $C.$ Note that $\Gamma^{*}\subset
\Gamma_E.$ Moreover, $\Gamma_{f(E)}>f(\Gamma^{*})$ and,
consequently,
\begin{equation}\label{eq7}
M\left(\Gamma_{f(E)}\right)\le M\left(f(\Gamma^{*})\right)\,.
\end{equation}
Consider $S_{\,\varepsilon}=S(x_0,\,\varepsilon),$ $
S_{\,\varepsilon_0}=S(x_0,\,\varepsilon_0),$
with $\varepsilon_0$ satisfying the assumption of lemma and
$\varepsilon\in\left(0,\,\varepsilon_0^{\,\prime}\right).$ Let
$R(r_1, r_2, x_0)=\{x\in {\Bbb R}^n: r_1<|x-x_0|<r_2\}.$ Note that
$\Gamma\left(S_{\varepsilon}, S_{\varepsilon_0}, R(\varepsilon,
\varepsilon_0, x_0)\right)<\Gamma^{\,*}$ and, consequently,
$f(\Gamma\left(S_{\varepsilon}, S_{\varepsilon_0}, R(\varepsilon,
\varepsilon_0, x_0)\right))<f(\Gamma^{*}).$ Hence
\begin{equation}\label{eq5}
M\left(f(\Gamma^{*})\right)\le
 M\left(f\left(\Gamma\left(S_{\varepsilon}, S_{\varepsilon_0},
 R(\varepsilon, \varepsilon_0, x_0)\right)\right)\right)\,.
\end{equation}
The relations (\ref{eq7}) and (\ref{eq5}) yield
$$
M\left(\Gamma_{f(E)}\right)\,\le
 M\left(f\left(\Gamma\left(S_{\varepsilon}, S_{\varepsilon_0},
 R(\varepsilon, \varepsilon_0, x_0)\right)\right)\right)
$$
and, consequently,
\begin{equation}\label{eq5aa}
{\rm cap}\,\,f(E) \le
 M\left(f\left(\Gamma\left(S_{\varepsilon}\,, S_{\varepsilon_0},
 R(\varepsilon, \varepsilon_0, x_0)\right)\right)\right)\,.
\end{equation}
Consider a family of Lebesgue measurable functions
$\eta_{\varepsilon}(t)=\psi_{\varepsilon}(t)/I(\varepsilon,
\varepsilon_0 ),$ $t\in(\varepsilon,\, \varepsilon_0).$ Given
$\varepsilon\in (0, \varepsilon_0^{\,\prime}),$ we have
$\int\limits_{\varepsilon}^{\varepsilon_0}\eta_{\varepsilon}(t)\,dt=1.$
By (\ref{eq2})
$$M\left(f\left(\Gamma\left(S_{\varepsilon}\,, S_{\varepsilon_0},
R(\varepsilon, \varepsilon_0, x_0)\right)\right)\right)\,\le$$
\begin{equation}\label{eq8}
\le\frac{1}{I^n(\varepsilon, \varepsilon_0)}
\int\limits_{\varepsilon<|x-x_0|<\varepsilon_0}Q(x)\cdot\psi_{\varepsilon}^n(|x-x_0|)\,
\ dm(x)
\end{equation}
for every $\varepsilon\in (0, \varepsilon_0^{\,\prime}).$ Now the
inequality (\ref{eq3B}) directly follows from (\ref{eq3.7B}),
(\ref{eq5aa}) and (\ref{eq8}). $\Box$

\medskip
A {\it chordal diameter} of a set $E\subset\overline{{\Bbb R}^n}$ is
the quantity
$$h(E)=\sup\limits_{x\,,y\in E}\,h(x,y)\,.$$
Denote by $\omega_{n-1}$ the area of the unit sphere ${\Bbb
S}^{n-1}$ in ${\Bbb R}^n.$

\begin{lemma}{}\label{lem6}{\,\sl
Let $f:D\rightarrow{\Bbb R}^n,$ $n\ge 2$ be an open discrete ring
$Q$-map\-ping at a point $x_0\in D,$ $D^{\,\prime}:=f(D)\subset B(0,
r)$ and $h\left(\overline{{\Bbb R}^n}\setminus B(0,
r)\right)\ge\delta>0.$ Suppose that there exist numbers $p \le n,$
$\varepsilon_0 \in (0,\, {\rm dist}\,(x_0,
\partial D)),$ $\varepsilon_0^{\,\prime}\in (0, \varepsilon_0)$ and nonnegative
Lebesgue measurable functions $\psi_{\varepsilon}:(\varepsilon,
\varepsilon_0)\rightarrow [0, \infty],$ $\varepsilon\in\left(0,
\varepsilon_0^{\,\prime}\right)$ such that
\begin{equation} \label{eq3.7A}
\int\limits_{\varepsilon<|x-x_0|<\varepsilon_0}Q(x)\cdot\psi_{\varepsilon}^n(|x-x_0|)
\ dm(x)\le K\cdot I^p(\varepsilon, \varepsilon_0)\qquad \forall\,\,
\varepsilon\in(0,\varepsilon_0^{\,\prime})\,,
\end{equation}
where the quantity $I(\varepsilon, \varepsilon_0)$ is defined by
(\ref{eq3}). Then
\begin{equation} \label{eq3.9}
h(f(x),f(x_0))\le\frac{\alpha_n}{\delta}\,\exp\{-\beta_n
I^{\gamma_{n,p}}(|x-x_0|, \varepsilon_0)\}
\end{equation}
for every $x\in B(x_0,{\varepsilon_0}^{\,\prime})$ with
%
\begin{equation}
\label{eq3.10A} \alpha_n=2\lambda_n^2,\quad \beta_n =
{\left(\frac{\omega_{n-1}}{K}\right)}^{\frac{1}{n-1}},\quad\gamma_{n,p}=1-\frac{p-1}{n-1},
\quad \lambda_n \in[4,2e^{n-1})\,.
\end{equation}
}
\end{lemma}

{\it Proof.} Choosing $\varepsilon\in (0, \varepsilon_0^{\,\prime})$
and $F(\varepsilon, \varepsilon_0):=K\cdot I^p(\varepsilon,
\varepsilon_0)$ one gets from (\ref{eq3B}) and (\ref{eq3.7A})
\begin{equation}\label{eq3a}
{\rm cap}\,f(E)\le K\cdot I^{p-n}\left(\varepsilon,
\varepsilon_0\right)\,.
\end{equation}
Applying \cite[Lemma~2.2]{Sev$_1$} and taking into account that
$f(A)\subset B(0, r),$ one obtains
\begin{equation}\label{eq17}
{\rm cap}\,f(E)\ge \frac {\omega_{n-1}}{\,\,\,\,{\left\{\log \frac
{2\lambda_n^2}{h \left(f(C)\right)h \left(\overline{{\Bbb
R}^n}\,\setminus B(0, r)\right)}\right\}}^{n-1}}
\end{equation}
where $\lambda_n \in [4, 2e^{n-1})$, $\lambda_2=4$ and
$\lambda_n^{1/n} \rightarrow e$ as $n \rightarrow \infty$ (cf.
\cite[(7.21) and Lemma~7.22]{Vu$_2$}). Combining (\ref{eq3a}) and
(\ref{eq17}) with assumption $h\left(\overline{{\Bbb R}^n}\setminus
B(0, r)\right)\ge\delta,$ we have
%
$$h\left(f(C)\right)\le\frac{2\,\lambda_n^2}{\delta}\,
\exp{\left\{\,-\,{\left(\frac{\omega_{n-1}}{K}\right)}^{\frac{1}{n-1}}
\,{\left(I(\varepsilon, \varepsilon_0
)\right)}^\frac{n-p}{n-1}\right\}}\,.$$
%
%
%
Letting $\alpha_n=2\lambda_n^2,\qquad \beta_n =
{\left(\frac{\omega_{n-1}}{K}\right)}^{\frac{1}{n-1}},\qquad\gamma_{n,p}=1-\frac{p-1}{n-1},$
\begin{equation}\label{eq20}
h\left(f(C)\right)\le\frac{\alpha_n}{\delta}\exp{\,\{-\beta_n
I^{\gamma_{n,p}}(\varepsilon, \varepsilon_0)\}}\,.
\end{equation}
%
Pick arbitrary $x\in D$ such that $|x-x_0|=\varepsilon,$ $0<
\varepsilon<\varepsilon_0^{\,\prime}.$ Then
$x\in\overline{B\left(x_0,\,\varepsilon\right)}$ and, moreover,
$f(x)\,\in\,f\left(\overline{B\left(x_0,\,\varepsilon\right)}\right)=f(C).$
Thus, the inequality holds for any $\varepsilon\in\,(0,
\varepsilon_0^{\,\prime}).$ Since
$\varepsilon\in\left(0,\varepsilon_0^{\,\prime}\right)$ is
arbitrary, the relation (\ref{eq3.9}) is fulfilled in the whole ball
$B(x_0, \varepsilon_0^{\,\prime}).$ $\Box$

\medskip
By the well known Liouville theorem, every conformal mapping of a
domain $D$ of ${\Bbb R}^n,$ $n\ge 3,$ is a restriction of a
M\"{o}bius transformations of a domain $D$ (see e.g. \cite{Hart} and
\cite[Theorem~2.5.I]{Ri}).  For $n=2,$ it is true when
$D:=\overline{{\Bbb R}^2}$ (see, e.g., \cite{Most}). Inversely, a
restriction of a M\"{o}bius transformation $U:\overline{{\Bbb
R}^n}\rightarrow \overline{{\Bbb R}^n}$ of a domain
$D\setminus\left\{\infty, U^{\,-1}(\infty)\right\}$ is conformal
mapping and, consequently, $M(U(\Gamma))=M(\Gamma)$ for every
M\"{o}bius transformation $U:\overline{{\Bbb R}^n}\rightarrow
\overline{{\Bbb R}^n}$ and any family of curves $\Gamma$ in
$\overline{{\Bbb R}^n}$ (cf. \cite[Theorem~8.1]{Va$_1$}).

\medskip
The following two lemmas describe the local behavior of the classes
of equicontinuous mappings obeying (\ref{eq2}).

\medskip
\begin{lemma}{}\label{lem1}{\,\sl Let ${\frak F}_Q$ be a family of open discrete
$Q$-map\-pings $f:D \rightarrow \overline{{\Bbb R}^n}$ at $x_0\in
D,$ $n \ge 2.$ Suppose that $p<n,$ $\varepsilon_0 \in (0,\, {\rm
dist}\,(x_0,
\partial D)),$ $\varepsilon_0^{\,\prime}\in (0, \varepsilon_0)$ and
$\psi:(0, \varepsilon_0)\rightarrow [0, \infty]$ is a nonnegative
Lebesgue measurable function  satisfying (\ref{eq3}) such that the
estimate (\ref{eq3.7A}) holds for $\psi_{\varepsilon}\equiv \psi$
and $I(\varepsilon, \varepsilon_0)\rightarrow \infty$ as
$\varepsilon\rightarrow 0.$

Then ${\frak F}_Q$ is equicontinuous at $x_0$ if and only if there
exist $\varepsilon_i=\varepsilon_i(x_0),$ $i=1, 2,$
$0<\varepsilon_1<\varepsilon_2<\varepsilon_0,$ such that the
inequality
\begin{equation}\label{eq22}
h(f(x),f(x_0))\le\alpha_n\exp\{-\widetilde{\beta_n}
I^{\gamma_{n,p}}(|x-x_0|, \varepsilon_2)\}\end{equation}
holds for every $f\in {\frak F}_Q$ and all $x\in B(x_0,
\varepsilon_1),$ where $\widetilde{\beta_n} =
{\left(\frac{\omega_{n-1}}{2K}\right)}^{\frac{1}{n-1}},$ $\alpha_n$
and $\gamma_{n,p}$ are defined by (\ref{eq3.10A}).}
\end{lemma}

\medskip
{\it Proof.} Since $I(\varepsilon, \varepsilon_0)\rightarrow \infty$
as $\varepsilon\rightarrow 0,$ the sufficient condition of this
lemma follows from Lemma \ref{lem6}. Now we show that inequality
(\ref{eq22}) holds when  ${\frak F}_Q$ is equicontinuous at $x_0.$
Given $\sigma>0,$ there exists $\Delta=\Delta(\sigma, x_0)$ such
that
\begin{equation}\label{eq4}
h(f(x), f(x_0))<\sigma\,,
\end{equation}
whenever $|x-x_0|<\Delta.$ Pick $\Delta<\varepsilon_0.$ Given a
mapping $f\in{\frak F}_Q$ and a point $x_0\in D,$ there is a
M\"{o}bius transformation $U:\overline{{\Bbb R}^n}\rightarrow
\overline{{\Bbb R}^n}$ such that
\begin{equation}\label{eq6}
U(f(x_0))=0, \qquad  h(U(f(x)), U(f(x_0)))=h (f(x), f(x_0))
\end{equation}
(see \cite[Theorem~12.2]{Va$_1$}). Clearly, the composition
$v:=U\circ f$ is ring $Q$-map\-ping at $x_0.$ For sufficiently small
$\sigma>0$ the relations (\ref{eq4}) and (\ref{eq6}) yield
$|v(x)|\le 1$ for some $\varepsilon_2=\varepsilon_2(x_0)$ and all
$x\in B(x_0, \varepsilon_2).$

Consider the mapping $g:=v|_{B(x_0, \varepsilon_2)}.$ Since
$I(\varepsilon, \varepsilon_0)\rightarrow \infty$ as
$\varepsilon\rightarrow 0,$  we have that
$0<I(\varepsilon, \varepsilon_0)<2I(\varepsilon,
\varepsilon_2)<\infty$
for all $\varepsilon\in (0, \varepsilon_1)$ and some
$\varepsilon_1\in (0, \varepsilon_2)$ and therefore by
(\ref{eq3.7A}) that
\begin{equation} \label{eq10}
\int\limits_{\varepsilon<|x-x_0|<\varepsilon_2}
Q(x)\cdot\psi^n(|x-x_0|) \ dm(x)\le 2K\cdot I^p(\varepsilon,
\varepsilon_2)\qquad \forall\,\, \varepsilon\in(0,\varepsilon_1)\,.
\end{equation}
Replacing $\varepsilon_2,$ $\varepsilon_1$ and $2K$ by
$\varepsilon_0,$ $\varepsilon_0^{\,\prime}$ and $K,$ we get the
inequality (\ref{eq3.7A}). Now applying Lemma \ref{lem6} to the
mapping $g$ with the conditions (\ref{eq6}) completes the proof.
$\Box$

\medskip
In what follows, we also need the following lemma.

\medskip
\begin{lemma}{}\label{lem2}{\,\sl Let ${\frak F}_Q$ be a family  of
open discrete mappings $f:D \rightarrow \overline{{\Bbb R}^n}$, $n
\ge 2,$ obeying (\ref{eq2}) at a point $x_0\in D.$ Suppose that
there exists a number $\widetilde{\varepsilon_0} \in (0,\, {\rm
dist}\,(x_0,
\partial D))$ such that the following property holds: for every
$\varepsilon_0 \in (0,\, \widetilde{\varepsilon_0})$ there exists
$\varepsilon_0^{\,\prime}\in (0, \varepsilon_0)$ and a nonnegative
Lebesgue measurable function $\psi:(0,
\widetilde{\varepsilon_0})\rightarrow [0, \infty]$ such that
\begin{equation}\label{eq13}
0<I(\varepsilon, \varepsilon_0):=
\int\limits_{\varepsilon}^{\varepsilon_0}\psi(t)dt <
\infty\qquad\forall\,\,\varepsilon\in(0,
\varepsilon_0^{\,\prime})\end{equation}
and
\begin{equation} \label{eq16}
\int\limits_{\varepsilon<|x-x_0|<\varepsilon_0}Q(x)\cdot\psi^n(|x-x_0|)
\ dm(x)\le K\cdot I^p(\varepsilon, \varepsilon_0)\qquad \forall\,\,
\varepsilon\in(0,\varepsilon_0^{\,\prime})\,.
\end{equation}
If the family ${\frak F}_Q$ is equicontinuous at $x_0,$ then there
are $\varepsilon_i=\varepsilon_i(x_0),$ $i=1, 2,$
$0<\varepsilon_1<\varepsilon_2<\varepsilon_0,$ such that
\begin{equation}\label{eq24}
h(f(x),f(x_0))\le\alpha_n\,\exp\{-\beta_n I^{\gamma_{n,p}}(|x-x_0|,
\varepsilon_2)\}
\end{equation}
for every $f\in {\frak F}_Q$ and all $x\in B(x_0, \varepsilon_1)$
with $\alpha_n,$ $\beta_n$ and $\gamma_{n,p}$ defined by
(\ref{eq3.10A}).}
\end{lemma}

\medskip
{\it Proof.} If the family ${\frak F}_Q$ is equicontinuous at the
point $x_0,$ then for every $\sigma>0$ there exists
$\Delta=\Delta(\sigma, x_0)$ such that (\ref{eq4}) holds for every
$x\in B(x_0, \Delta).$ Given a mapping $f\in {\frak F}_Q,$ there is
a M\"{o}bius transformation $U:\overline{{\Bbb R}^n}\rightarrow
\overline{{\Bbb R}^n}$ provided (\ref{eq6}). The mapping $v:=U\circ
f$ is a ring $Q$-map\-ping at $x_0.$ For sufficiently small
$\sigma>0$ the relations (\ref{eq4}) and (\ref{eq6}) yield
$|v(x)|\le 1$ for small enough $\sigma,$ for some
$\varepsilon_2=\varepsilon_2(x_0)$ and all $x\in B(x_0,
\varepsilon_2).$ Consider the mapping $g:=f|_{B(x_0,
\varepsilon_2)}.$ By (\ref{eq13}) the double inequality
$0<I(\varepsilon, \varepsilon_2)<\infty$ holds for all
$\varepsilon\in (0, \varepsilon_1)$ and some $\varepsilon_1\in (0,
\varepsilon_2)$ and therefore by (\ref{eq16})
\begin{equation} \label{eq23}
\int\limits_{\varepsilon<|x-x_0|<\varepsilon_2}
Q(x)\cdot\psi^n(|x-x_0|) \ dm(x)\le K\cdot I^p(\varepsilon,
\varepsilon_2)\qquad \forall\,\, \varepsilon\in(0,\varepsilon_1)\,.
\end{equation}
Again applying Lemma \ref{lem6} to the mapping $g$ with the
conditions (\ref{eq6}) completes the proof. $\Box$

\medskip
\begin{remark}{}\label{remark1}
The inequalities (\ref{eq22}) and (\ref{eq24}) in Lemmas \ref{lem1}
and \ref{lem2} are both of the same type and provide the upper
bounds for the chordal distances between $f(x)$ and $f(x_0).$ The
constant $\beta_n=(\frac{\omega_{n-1}}{K})^{1/(n-1)}$ in
(\ref{eq24}) is better than the corresponding constant
$\widetilde{\beta_n}=(\frac{\omega_{n-1}}{2K})^{1/(n-1)}$ in
(\ref{eq22}). On the other hand, the conditions (\ref{eq3}),
(\ref{eq3.7A}) of Lemma \ref{lem1} hold for one fixed
$\varepsilon_0,$ while the assertion of Lemma \ref{lem2} is
fulfilled for infinitely many $\varepsilon_0$ with the same $K>0.$
Applications of both lemmas will be presented in the next section.
\end{remark}

\medskip
\section{On necessary and sufficient conditions of
equicontinuity} \setcounter{equation}{0}

\medskip
We start to prove Theorem \ref{th3} by applying Lemma \ref{lem1} for
an appropriate function $\psi.$ Letting $\psi(t):=\frac{1}{t\log
1/t}$ and $\varepsilon_0<1/e,$ we have that
\begin{equation}\label{eq32}I(\varepsilon,
\varepsilon_0)=\int\limits_{\varepsilon}^{\varepsilon_0}\psi(t)\,dt=
\log{\frac{\log{\frac{1}
{\varepsilon}}}{\log{\frac{1}{\varepsilon_0}}}}\,,
\end{equation}
where $I(\varepsilon, \varepsilon_0)\rightarrow \infty$ as
$\varepsilon\rightarrow 0$ and the condition (\ref{eq3.7A}) holds by
\cite[Corollary~2.3]{IR} (cf. \cite[Lemma~6.1]{MRSY}). Thus, the
condition (\ref{eq1}) follows from (\ref{eq3.9}) for the chosen
function $\psi.$ $\Box.$

\medskip
Given a Lebesgue measurable function $Q:D\rightarrow [0,\infty],$
denote by $q_{x_0}(r)$ the integral average of $Q(x)$ over the
sphere $|x-x_0|=r,$
\begin{equation}\label{eq3.1}
q_{x_0}(r):=\frac{1}{\omega_{n-1}r^{n-1}}\int\limits_{|x-x_0|=r}Q(x)\,dS\,,
\end{equation}
where $dS$ denotes an element of the sphere area. The following
statements are relied on Lemma \ref{lem2}.

\medskip
\begin{theorem}{}\label{th4}{\,\sl Given $x_0\in D$ and a Lebesgue measurable function
$Q:D\rightarrow [0, \infty],$ suppose that the condition
(\ref{eq12}) holds for some $\varepsilon_0=\varepsilon_0(x_0)\in (0,
{\rm dist\,}(x_0,
\partial D)).$ A family of open
discrete ring $Q$-map\-pings $f:D\rightarrow \overline{{\Bbb R}^n}$
at the point $x_0\in D$ is equicontinuous at $x_0$ if and only if
the inequality
\begin{equation}\label{eq11}
h(f(x), f(x_0))\le {\alpha}_n\cdot
\exp\left\{-\int\limits_{|x-x_0|}^{\varepsilon_0}
\frac{dt}{tq_{x_0}^{\,\frac{1}{n-1}}(t)}\right\}
\end{equation}
holds for all $x\in B(x_0, \varepsilon_0^{\,\prime})$ and some
$\varepsilon_0^{\,\prime}\in (0, \varepsilon_0)$ with $\alpha_n$
defined in Lemma \ref{lem6}. }
\end{theorem}

\medskip
{\it Proof.} The sufficient condition of Theorem \ref{th4} is
established by a direct application of the inequality (\ref{eq11})
and the condition (\ref{eq12}) from Theorem \ref{th1}. So, it
remains to prove the necessity. To this end, pick
\begin{equation}\label{eq36} \psi(t)\quad=\quad \left
\{\begin{array}{rr} 1/[tq^{\frac{1}{n-1}}_{x_0}(t)]\ , & \ t\in
(\varepsilon, \varepsilon_0),
\\ 0\ ,  &  \ t\notin (\varepsilon,
\varepsilon_0)\,.
\end{array} \right.
\end{equation}
By (\ref{eq12}), for sufficiently small $\varepsilon_0>0$ there
exists $\varepsilon_0^{\,\prime}\in (0, \varepsilon_0)$ such that
$I(\varepsilon,
\varepsilon_0):=\int\limits_{\varepsilon}^{\varepsilon_0}\psi(t)dt>0$
for every $\varepsilon\in (0, \varepsilon_0^{\,\prime}),$ and,
moreover, $I(\varepsilon, \varepsilon_0)<\infty$ for every
$\varepsilon\in (0, \varepsilon_0)$ (see \cite[Remark~3]{Sev$_2$}).
A direct calculation yields
\begin{equation}\label{eq37}\int\limits_{\varepsilon<|x-x_0|<\varepsilon_0}
Q(x)\cdot\psi^n(|x-x_0|)\ dm(x)=\omega_{n-1}\cdot I(\varepsilon,
\varepsilon_0)\,.
\end{equation}
Thus, $\psi$ satisfies condition (\ref{eq16}) with $p=1$ and
$K=\omega_{n-1},$ and the desired conclusion follows from Lemma
\ref{lem2}. $\Box$

\medskip
The following statement follows directly from Theorem \ref{th4}.

\medskip
\begin{corollary}{}\label{cor1}
{\,\sl Given $x_0\in D$ and a Lebesgue measurable function
$Q:D\rightarrow [0, \infty],$ assume that the estimate
\begin{equation}\label{eq40} q_{x_0}(r)\le C
\left(\log\frac{1}{r}\right)^{n-1}\qquad\forall\,\,r\in (0,
\varepsilon_0)
\end{equation}
holds for some $\varepsilon_0=\varepsilon(x_0)\in (0, {\rm
dist\,}(x_0,
\partial D)).$ Then a family of all open discrete ring $Q$-map\-pings $f:D\rightarrow
\overline{{\Bbb R}^n}$ at $x_0$  is equicontinuous at $x_0$ if and
only if
\begin{equation}\label{eq39}
h(f(x), f(x_0))\le
\frac{M}{\log^{\left(\frac{1}{C}\right)^{\frac{1}{n-1}}}\frac{1}{|x-x_0|}}\quad\forall\,\,
x\in B(x_0, \varepsilon_0^{\,\prime})
\end{equation}
and for some $\varepsilon_0^{\,\prime}\in (0, \varepsilon_0),$ where
the constant $M$ depends only on $n$ and $x_0.$}
\end{corollary}

\medskip
The next statement shows that the condition (\ref{eq12}) in Theorem
(\ref{th4}) cannot be removed.

\medskip
\begin{theorem}{}\label{th5}
{\sl\, Let $D$ be a domain in ${\Bbb R}^n,$ $n\ge 2,$ $x_0\in D,$
and $0<\varepsilon_0<{\rm dist\,}(x_0, \partial D).$ Given a
function $Q:D\rightarrow [1, \infty],$ $Q\in L_{loc}^1(D),$ with
$$\int\limits_{0}^{\varepsilon_0}\frac{dt}{tq_{x_0}^{\,\frac{1}{n-1}}(t)}<\infty\,,$$
where $q_{x_0}$ is defined in (\ref{eq3.1}). Then there exists a
family of uniformly bounded ring $Q$-map\-pings at the point $x_0$
which is not equicontinuous at $x_0.$}
\end{theorem}

{\it Proof.} Without loss of generality, we may assume that $D={\Bbb
B}^n=\{x\in{\Bbb R}^n: |x|<1\}$ and $x_0=0.$ Define a sequence of
mappings $f_m:{\Bbb B}^n\rightarrow {\Bbb R}^n$ by
$$f_m(x)=\frac{x}{|x|}\,\rho_m(|x|)\,,\qquad f_m(0):=0\,,$$
where
$$\rho_m(r)=
\exp\left\{-\int\limits_{r}^1\frac{dt}{tq_{0,
m}^{1/(n-1)}(t)}\right\}\,, \qquad q_{0,
m}(r):=\frac{1}{\omega_{n-1}r^{n-1}}\int\limits_{|x|=r}Q_m(x)\,dS\,,
$$
$$Q_m(x)\quad=\quad \left \{\begin{array}{rr} Q(x) , & \
|x|> 1/m\ ,
\\ 1\ ,  &  |x|\le 1/m\,.
\end{array} \right.$$
Now we show that every $f_m,$ $m=1,2,\ldots,$ is a ring
$Q$-ho\-me\-o\-mor\-phism at the point $x_0=0.$ Obviously, $f(S(0,
r))=S(0, R_m)$ where
$$R_m=\exp\left\{-\int\limits_{r}^1\frac{dt}{tq_{0,
m}^{1/(n-1)}(t)}\right\}\,,\quad r\in (0, 1)\,.$$ Note that
$$f_m(\Gamma(S(0, r_1), S(0,
r_2), R(r_1, r_2, 0)))= \Gamma(S(0, R_{1, m}), S(0, R_{2, m}),
R(R_{1, m}, R_{2, m}, 0))\,,$$ where
$R_{i, m}=\exp\left\{-\int\limits_{r_i}^1\frac{dt}{tq_{0,
m}^{1/(n-1)}(t)}\right\},$ $i=1, 2.$ Following to
\cite[7.5]{Va$_1$},
$$M(f(\Gamma(S(0, r_1), S(0, r_2), R(r_1, r_2, 0))))=\frac{\omega_{n-1}}
{\left(\int\limits_{r_1}^{r_2}\frac{dt}{tq_{0,
m}^{1/(n-1)}(t)}\right)^{n-1}}\le \frac{\omega_{n-1}}
{\left(\int\limits_{r_1}^{r_2}\frac{dt}{tq_{0}^{1/(n-1)}(t)}\right)^{n-1}}\,,$$
where $q_0(r)$ is defined by (\ref{eq3.1}) at the origin.
Consequently, the mapping $f_m$ is a ring $Q$-map\-ping at $0$ (see
\cite[Theorem~1]{Sev$_3$}). Since $|f_m(x)|\le 1$ for every $m\in
{\Bbb N},$ the family $\{f_m(x)\}_{m=1}^{\infty}$ is uniformly
bounded. On the other hand, for every sequence $x_m$ provided
$|x_m|=1/m,$ $m=1,2,\ldots,$ one obtains $|f_m(x_m)|\ge \sigma,$
where $\sigma>0$ does not depend on $m.$ Thus, the family
$\{f_m(x)\}_{m=1}^{\infty}$ is not equicontinuous at the origin.
$\Box$

\medskip
\section{On equicontinuity of map families omitting sets }
\setcounter{equation}{0}

The main results of this section are relied on Vuorinen's approach
in \cite{Vu$_1$}. Here we prove more general results on
equicontinuity of open discrete mappings satisfying the estimates
(\ref{eq2})--(\ref{eq2A}) (cf. \cite{Sev$_1$}). Following
\cite{Vu$_1$}, we establish these results on equicontinuity of the
mapping families involving set function of a special kind.

\medskip
Let $B^{\,*}(x, t)=\{y\in \overline{{\Bbb R}^n}: h(x, y)<t\}$ be a
spherical ball centered at the point $x$ of the radius $t.$ Given
$x\in \overline{{\Bbb R}^n},$ $E\subset \overline{{\Bbb R}^n}$ and
$0<r<t<1, $ we define
$$\left \{\begin{array}{rr}
m_t(E, r, x)= M(\Gamma(\partial B^{\,*}(x, t), E\cap
\overline{B^{\,*}(x, r)}))\,,
\\ m(E, x)=m_{\sqrt{3}/2}(E, \frac{\sqrt{2}}{2}, x)\,,
\end{array} \right.$$
and
\begin{equation}\label{eq25}
\left \{\begin{array}{rr} c(E, x)=\max\{m(E, x), m(E,
\widetilde{x})\}\,,
\\ c(E)=\inf\limits_{x\in \overline{{\Bbb R}^n}}c(E, x)\,,
\end{array} \right.
\end{equation}
where $\widetilde{x}=-x/|x|^2.$ Given a compact set
$E\subset\overline{{\Bbb R}^n},$ we have $c(E)=0$ if and only if
${\rm cap\,}E=0$ (see \cite[Corollary~3.19]{Vu$_1$}).

\medskip
Let $D$ be a domain in ${\Bbb R}^n,$ $n\ge 2.$ Denote by ${\frak
F}_{Q, \delta}$ a family of all open discrete ring $Q$-map\-pings
$f:D\rightarrow\overline{{\Bbb R}^n}\setminus E_f$ at $x_0\in D$
with $c(E_f)\ge \delta>0,$ where $E_f$ is compact and $c(\cdot)$ is
defined by (\ref{eq25}).

\medskip
The following auxiliary result is of independent interest.

\medskip
\begin{lemma}{}\label{lem9} {\sl\, Suppose that there exist
 $\varepsilon_0>0,$ $\varepsilon_0<{\rm dist\,}(x_0,
\partial D),$ $\varepsilon_1\in (0, \varepsilon_0)$
and a function $\psi:(0, \varepsilon_0)\rightarrow [0, \infty]$ with
\begin{equation}\label{eq38}0<I(\varepsilon,
\varepsilon_0)=\int\limits_{\varepsilon}^{\varepsilon_0} \psi(t)dt <
\infty\qquad\forall\,\,\varepsilon \in(0, \varepsilon_1)
\end{equation}
such that
\begin{equation} \label{eq26}
\int\limits_{\varepsilon<|x-x_0|<\varepsilon_0}
Q(x)\cdot\psi^n(|x-x_0|) \ dm(x)=o\left(I^n\left(\varepsilon,
\varepsilon_0\right)\right)
\end{equation}
%
%
as $\varepsilon\rightarrow 0.$ Then the family ${\frak F}_{Q,
\delta}$ is equicontinuous at $x_0.$}
\end{lemma}

\medskip
{\it Proof.} Given a mapping $f\in {\frak F}_{Q, \delta},$ consider
the condensers $E=(A, C)$ and $E^{\,\prime}=f(E)=(f(A), f(C))$ where
$C:=\overline{B(x_0, \varepsilon)},$ $\varepsilon \in(0,
\varepsilon_1),$ $A=B(x_0, r_0),$ $r_0={\rm dist\,}(x_0,
\partial D).$ Following the notation $\Gamma_E$ in the proof
of Lemma \ref{lem4}, let $\Gamma_{E^{\,\prime}}$ be a corresponding
curve family for the condenser $E^{\,\prime}.$ Note that by
\cite[Theorem~1.I.46]{Ku$_2$} $\Gamma(f(C),
E_f)>\Gamma_{E^{\,\prime}}$ and, consequently, by
\cite[Theorem~6.4]{Va$_1$}), Lemma \ref{lem4} and (\ref{eq26}), one
gets
\begin{equation}\label{eq27}
M(\Gamma(f(C), E_f))\le M(\Gamma_{E^{\,\prime}})={\rm cap\,}f(E)\le
\alpha(\varepsilon)\,,
\end{equation}
where $\alpha(\varepsilon)\rightarrow 0$ as $\varepsilon\rightarrow
0.$ On the other hand, by \cite[Theorem~3.14]{Vu$_1$},
\begin{equation}\label{eq28}
M(\Gamma(f(C), E_f))\ge\beta_n\min\{c(f(C)), c(E_f)\}\,,
\end{equation}
where the constant $\beta_n$ depends only on $n.$

Since $c(F)\ge a_n h(F)$ for every connected set $F$ in
$\overline{{\Bbb R}^n},$
\begin{equation}\label{eq29}
c(f(C))\ge a_n \cdot h(f(C))
\end{equation}
where $h(F)$ is the chordal diameter of $F$ and $a_n$ is some
constant (see \cite[Corollary~3.13]{Vu$_1$}). It is known that
$c(E)\le \omega_{n-1}\cdot\left(\log\sqrt{3}\right)^{1-n}$ for every
set $E\subset \overline{{\Bbb R}^n},$ consequently,
\begin{equation}\label{eq30}
\frac{c(E)}{\omega_{n-1}\cdot\left(\log\sqrt{3}\right)^{1-n}}\quad\le\quad
1\,,\end{equation}
cf. \cite[(3.7)]{Vu$_1$} and (\ref{eq25}).
If the minimum in (\ref{eq28}) equals $c(f(C)),$ then we have by
(\ref{eq29}) and (\ref{eq30})
\begin{equation}\label{eq31}
M(\Gamma(f(C), E_f))\ge\beta_n\cdot a_n \cdot h(f(C))\ge
\frac{\beta_n\cdot a_n \cdot
h(f(C))c(E_f)}{\omega_{n-1}\cdot\left(\log\sqrt{3}\right)^{1-n}}\,.
\end{equation}
For the case when $\min\{c(f(C)), c(E_f)\}=c(E_f),$ the inequality
(\ref{eq28}) yields
\begin{equation}\label{eq32A} M(\Gamma(f(C), E_f))\ge \beta_n\cdot c(E_f)\ge
\beta_n\cdot h(f(C))c(E_f)\,.\end{equation}
Letting $c_n:=\min\left\{\beta_n, \frac{\beta_n\cdot a_n
}{\omega_{n-1}\cdot\left(\log\sqrt{3}\right)^{1-n}}\right\},$ one
can derive from (\ref{eq31}) and (\ref{eq32A}) that
\begin{equation}\label{eq33}
M(\Gamma(f(C), E_f))\ge c_n\cdot h(f(C))c(E_f)\ge
c_n\cdot\delta\cdot h(f(C))
\end{equation}
and combining (\ref{eq27}) and (\ref{eq33}) the following estimate
holds
\begin{equation}\label{eq35}
h(f(C))\le \frac{\alpha(\varepsilon)}{c_n\delta}\,.
\end{equation}
Since $\alpha(\varepsilon)\rightarrow 0$ as $\varepsilon\rightarrow
0,$ one can conclude from (\ref{eq35}) that for every $\sigma>0$
there exists $\Delta=\Delta(\sigma)$ such that $h(f(C))<\sigma$
whenever $\varepsilon<\Delta.$ In other words, $h(f(x),
f(x_0))<\sigma$ for all $x\in B(x_0, \varepsilon),$
$\varepsilon<\Delta.$ Hence, the inequality $h(f(x), f(x_0))<\sigma$
holds for all $x\in B(x_0, \Delta)$ and every $f\in {\frak F}_{Q,
\Delta},$ i.e. the family ${\frak F}_{Q, \Delta}$ is equicontinuous
at the point $x_0.$ $\Box$

\medskip
The following result is a strengthen version of
 \cite[Lemma~8]{Sev$_2$}.

\medskip
\begin{lemma}{}\label{pr1}
{\,\sl Let $Q:D\rightarrow [0, \infty]$ be a Lebesgue measurable
function. Assume that either of the following conditions holds at
$x_0\in D:$ 1) $Q\in FMO(x_0);$ 2) $q_{x_0}(r)\le C\cdot
(\log\frac{1}{r})^{n-1};$ 3) the relation (\ref{eq12}) is true for
some $\varepsilon_0=\varepsilon_0(x_0).$ Then the conditions
(\ref{eq38})--(\ref{eq26}) of Lemma \ref{lem9} hold.}
\end{lemma}

\medskip
{\it Proof.} We start with the assumption 1) choosing
$\psi(t):=\frac{1}{t\log 1/t}$ and $\varepsilon_0<e^{\,-1}.$ Then
the conditions (\ref{eq32}) and (\ref{eq26}) hold (see
\cite[Lemma~6.1]{MRSY}). Thus, the relations
(\ref{eq38})--(\ref{eq26}) are fulfilled.

Since the assumption 2) is a particular case of 3), we should prove
the implication (\ref{eq12})$\Rightarrow$
(\ref{eq38})--(\ref{eq26}). Given sufficiently small
$\varepsilon>0,$ pick $\psi$ by (\ref{eq36}). In accordance to
(\ref{eq12}), there exists $\varepsilon_1\in (0, \varepsilon_0)$
provided $I(\varepsilon, \varepsilon_0)>0$ for all $\varepsilon\in
(0, \varepsilon_1).$ Note that $I(\varepsilon,
\varepsilon_0)<\infty$ for every $\varepsilon\in (0, \varepsilon_0)$
(see \cite[Remark~3]{Sev$_2$}). Now arguing similar to the proof of
Theorem \ref{th4}, one obtains equality (\ref{eq37}) and, moreover,
the asymptotic relation (\ref{eq26}). $\Box$

\medskip
The assertion of Theorem \ref{th1} directly follows from Lemma
\ref{lem9} and Lemma \ref{pr1}.

\section{Applications and related results}
\setcounter{equation}{0}

Now we present some important particular cases of Theorem \ref{th1}.
The first one relates to the case $Q(x)$ is uniformly bounded by a
constant $K,$ i.e., $Q(x)\le K=const$ and $E_f=E\subset
\overline{{\Bbb R}^n}.$ Thus, we recall the well--known results on
equicontinuity for quasiregular mappings (see
\cite[Theorem~3.17]{MRV$_2$} and \cite[Corollary~2.7.III]{Ri}).
\medskip

\begin{corollary}{}\label{cor2}
{\sl\, Given a compact set $E\subset \overline{{\Bbb R}^n}$ with
${\rm cap\,}E>0,$ the family of all $K$--qua\-si\-re\-gular mappings
$f:D\rightarrow\overline{{\Bbb R}^n}\setminus E$ is equicontinuous
at every point $x_0\in D.$}
\end{corollary}

\medskip
As a consequence for $E_f\equiv E\subset \overline{{\Bbb R}^n},$ we
obtain the following author result (see
\cite[Theorems~5.1--5.2]{Sev$_1$}).
\medskip

\begin{corollary}{}\label{cor3}
{\sl\, Given a compact set $E\subset \overline{{\Bbb R}^n}$ with
${\rm cap\,}E>0,$ a family of all open discrete $Q$-map\-pings
$f:D\rightarrow\overline{{\Bbb R}^n}\setminus E$ is equicontinuous
at every point $x_0\in D$ if either of the conditions 1)--3) of
Theorem \ref{th1} holds.}
\end{corollary}

\medskip
The following result is based on Theorems \ref{th3}, \ref{th1},
\ref{th4} and Corollary \ref{th4}.

\medskip
\begin{theorem}{}\label{th6}{\sl\, Let $D$ be a domain in ${\Bbb R}^n,$ $n\ge 2.$ Denote by ${\frak
F}_{Q, \delta}$ a family of all open discrete ring $Q$-map\-pings
$f:D\rightarrow\overline{{\Bbb R}^n}\setminus E_f$ at $x_0\in D$
with $c(E_f)\ge \delta>0,$ where $E_f$ is compact and $c(\cdot)$ is
defined by (\ref{eq25}). Then:

 1. Every $f\in {\frak F}_{Q, \delta}$ satisfies the estimate
(\ref{eq1}) for all $x\in B(x_0, \varepsilon_0(x_0))$ and some
$\varepsilon_0(x_0)\in (0, {\rm dist\,}(x_0, \partial D))$ whenever
$Q\in FMO(x_0)$ where $p=p(n , Q)>0$ and $C_n>0$ are some constants.
 и всех

\medskip
2. Every $f\in {\frak F}_{Q, \delta}$ satisfies the estimate
(\ref{eq11}) for some constant $\alpha_n>0$ whenever the condition
(\ref{eq12}) holds for some $\varepsilon_0\in (0, {\rm dist\,}(x_0,
\partial D))$

\medskip

3. Every $f\in {\frak F}_{Q, \delta}$ satisfies the estimate
(\ref{eq39}) in a neighborhood of the point $x_0$ whenever the
inequality (\ref{eq40}) holds with $C>0$ where $M>0$ is some
constant depending only on $n$ and point $x_0.$}
\end{theorem}

\medskip
Finally, we give important applications of our main results to the
Sobolev classes. Let $x\in D$ be a point of differentiability of
$f.$ Set
$$l\left(f^{\,\prime}(x)\right)\,=\,\min\limits_{h\in {\Bbb
R}^n \backslash \{0\}} \frac {|f^{\,\prime}(x)h|}{|h|}\,, \Vert
f^{\,\prime}(x)\Vert\,=\,\max\limits_{h\in {\Bbb R}^n \backslash
\{0\}} \frac {|f^{\,\prime}(x)h|}{|h|}\,, J(x,f)={\rm det}
f^{\,\prime}(x),$$
and define for such $x\in D$ the {\it inner dilatation of $f$ at
$x$} by
$$K_{I}(x,f)\quad =\quad\left\{
\begin{array}{rr}
\frac{|J(x,f)|}{{l\left(f^{\,\prime}(x)\right)}^n}, & J(x,f)\ne 0,\\
1,  &  f^{\,\prime}(x)=0, \\
\infty, & {\rm otherwise}
\end{array}
\right.\,.$$

\medskip Recall that a point $y_0\in D$ is said to be a {\it branch
point} of a mapping $f:D\rightarrow {\Bbb R}^n$ if, for every
neighborhood $U$ of $y_0,$ a restriction $f|_{U}$ fails to be
homeomorphic. A set of all branch points of $f$ is denoted by $B_f.$

\medskip
The following statement follows from \cite[Corollary 2]{Sev$_4$} and
Theorems \ref{th3}, \ref{th4} and Corollary \ref{cor1} as well. Let
$\frak{A}_{Q}$ be a family of all discrete open mappings
$f:D\rightarrow{\Bbb R}^n$ of Sobolev's class $W_{loc}^{1, n}$ with
a branch set $B_f$ of the Lebesgue measure zero. Assume that $K_I(x,
f)\le Q(x)\in L_{loc}^1.$ Then we have:

\medskip
\begin{corollary}{}\label{cor4}
{\sl\, 1. Every $f\in\frak{A}_{Q}$ satisfies the estimate
(\ref{eq1}) for all $x\in B(x_0, \varepsilon_0(x_0))$ and some
$\varepsilon_0(x_0)\in (0, {\rm dist\,}(x_0, \partial D))$ whenever
$Q\in FMO(x_0)$ where $p=p(n , Q)>0$ and $C_n>0$ are some constants.

\medskip
2. Every $f\in \frak{A}_{Q}$ satisfies the estimate (\ref{eq11})
with some constant $\alpha_n>0$ whenever the condition (\ref{eq12})
holds for some $\varepsilon_0\in (0, {\rm dist\,}(x_0, \partial
D)).$

\medskip
3. Every $f\in\frak{A}_{Q}$ satisfies the estimate (\ref{eq39}) in a
neighborhood of the point $x_0$ whenever the inequality (\ref{eq40})
holds with $C>0$ where $M>0$ is some constant depending only on $n$
and point $x_0.$}
\end{corollary}

\medskip
Denote by $\frak{B}_{Q, \delta}$ a family of all discrete open
mappings $f:D\rightarrow{\Bbb R}^n\setminus E_f$ of $W_{loc}^{1, n}$
with a branch set $B_f$ of the Lebesgue measure zero such that
$K_I(x, f)\le Q(x)\in L_{loc}^1$ and $c(E_f)\ge\delta>0.$ The
following statement follows from \cite[Corollary~2]{Sev$_4$} and
Theorems~\ref{th3}, \ref{th4}, and Corollary \ref{cor1}.

\medskip
\begin{corollary}{}\label{cor5}
{\sl\, A family of the mappings $\frak{B}_{Q, \delta}$ is
equicontinuous at $x_0$ provided that either of the following
conditions  1)--3) of Theorem \ref{th1} holds.}
\end{corollary}

\medskip
The modulus and capacity inequalities provide one of the main tools
in investigation the basic properties of a wide spectrum of mappings
(see, e.g., \cite{AC}--\cite{GS}, \cite{IR}, \cite{LV},
\cite{MRV$_1$}--\cite{M}, \cite{Pol}, \cite{Ri}--\cite{Vu$_2$}).
Thus, the results of the paper can be successfully applied to these
classes of mappings.

\medskip
\noindent
{\bf Evgeny Sevost'yanov:}\\
Institute of Applied Mathematics and Mechanics,\\
National Academy of Sciences of Ukraine,\\ 74 Roze Luxemburg Str.,
Donetsk,\\ 83114, UKRAINE\\
e--mail: brusin2006@rambler.ru

\end{document}